\documentclass[12pt]{article}
\usepackage{amssymb}
\usepackage{amsfonts}
\usepackage{graphicx}
\usepackage{amsmath}
\usepackage{cite}
\usepackage{ifthen}
\usepackage{float}

\newtheorem{theorem}{Theorem}
\newtheorem{algorithm}[theorem]{Algorithm}

\newtheorem{proposition}[theorem]{Proposition}
\newtheorem{remark}[theorem]{Remark}

\newenvironment{proof}[1][{\bf Proof}]{\noindent \textbf{#1.} }{\ \rule{0.5em}{0.5em}}

\floatstyle{ruled}
\newfloat{algorithm}{htb}{alg}
\floatname{algorithm}{Algorithm}
\newcounter{algo@row}
\newcounter{algo@rowindent}
\newcommand{\algofont}[1]{\textbf{#1}}
\newcommand{\algonumbersize}[1]{\scriptsize{#1}}
\newcommand{\algopreitem}[1][\arabic{algo@row}]{\texttt{\algonumbersize{#1}}}
\newcommand{\algoitemskip}{\hspace{\value{algo@rowindent}cc}}

\newenvironment{algo}{\vskip.3em\small \begin{list}{\algopreitem\texttt{\algonumbersize{:}}}{      \usecounter{algo@row}      \setcounter{algo@rowindent}{0}      \setlength{\itemindent}{2em}      \setlength{\labelwidth}{2em}      \setlength{\parsep}{0cm}    }}{
  \end{list}\vskip-.5em
}
\newcommand{\algonewnestedopen}[2]{
  \newcommand{#1}[1][]{    \ifthenelse{\equal{##1}{}}{\item}{\item[{\algopreitem[##1]}]}
    \algoitemskip\algofont{#2}    \addtocounter{algo@rowindent}{1}    \ignorespaces
  }
}
\newcommand{\algonewnestedaux}[2]{
  \newcommand{#1}[1][]{
    \addtocounter{algo@rowindent}{-1}
    \ifthenelse{\equal{##1}{}}{\item}{\item[{\algopreitem[##1]}]}
    \algoitemskip\algofont{#2}    \addtocounter{algo@rowindent}{+1}    \ignorespaces
  }
}
\newcommand{\algonewnestedclose}[2]{
  \newcommand{#1}[1][]{
    \addtocounter{algo@rowindent}{-1}
    \ifthenelse{\equal{##1}{}}{\item}{\item[{\algopreitem[##1]}]}
    \algoitemskip\algofont{#2}    \ignorespaces
  }
}
\newcommand{\algonewcommand}[2]{
  \newcommand{#1}[1][default]{
    \ifthenelse{\equal{##1}{default}}{\item}{\item[{\algopreitem[##1]}]}    \algoitemskip\algofont{#2}    \ignorespaces
  }}
\newcommand{\algonewkeyword}[2]{\newcommand{#1}{\algofont{#2}}}
\algonewcommand{\STATE}{\ignorespaces}
\algonewcommand{\REQUIRE}{Require }
\algonewcommand{\DEFINE}{Define }
\algonewcommand{\COMPUTE}{Compute }
\algonewcommand{\UPDATE}{Update }
\algonewcommand{\OUTPUT}{Output }
\algonewcommand{\COMMENT}{\hfill // }
\algonewnestedopen{\IF}{if }
\algonewnestedaux{\ELSEIF}{else if }
\algonewnestedaux{\THEN}{then }
\algonewnestedaux{\ELSE}{else}
\algonewnestedclose{\ENDIF}{end if }
\algonewnestedopen{\FOR}{for }
\algonewnestedclose{\ENDFOR}{end for }
\algonewkeyword{\To}{to }
\algonewkeyword{\Do}{do }
\algonewkeyword{\If}{if }
\algonewkeyword{\Then}{then }
\algonewkeyword{\Else}{else }
\algonewkeyword{\Endif}{end if }
\begin{document}

\title{A rational Arnoldi approach for ill-conditioned linear systems}
\author{C. Brezinski\thanks{%
Laboratoire Paul Painlev\'e, UMR CNRS 8524, UFR de Math\'ematiques Pures et
Appliqu\'ees, Universit\'e des Sciences et Technologies de Lille,
59655--Villeneuve d'Ascq cedex, France. E--mail: \texttt{%
Claude.Brezinski@univ-lille1.fr}} \and P. Novati \thanks{
Universit\`a degli Studi di Padova, Dipartimento di Matematica Pura ed
Applicata, Via Trieste 63, 35121--Padova, Italy. E--mail: \texttt{%
novati@math.unipd.it}} \and M. Redivo--Zaglia \thanks{
Universit\`a degli Studi di Padova, Dipartimento di Matematica Pura ed
Applicata, Via Trieste 63, 35121--Padova, Italy. E--mail: \texttt{%
Michela.RedivoZaglia@unipd.it}} }
\maketitle

\begin{abstract}
For the solution of full-rank ill-posed linear systems a new approach based
on the Arnoldi algorithm is presented. Working with regularized systems, the
method theoretically reconstructs the true solution by means of the
computation of a suitable function of matrix. In this sense the method can
be referred to as an iterative refinement process. Numerical experiments
arising from integral equations and interpolation theory are presented.
Finally, the method is extended to work in connection with the standard
Tikhonov regularization with a right hand side contaminated by noise.
\end{abstract}


\section{Introduction}

In this paper we consider the solution of ill-conditioned linear systems%
\begin{equation}
Ax=b.  \label{pr}
\end{equation}%
We mainly focus the attention on linear systems in which $A\in \mathbb{R}%
^{N\times N}$ is full rank with singular values that gradually decay to $0$,
as for instance in the case of the discretized Fredholm integral equations
of the first kind. In order face this kind of problems one typically apply
some regularization technique such as the well known Tikhonov regularization
(see e.g. \cite{PCH} for a wide background). The Tikhonov regularized system
takes the form%
\begin{equation}
(A^{T}A+\lambda H^{T}H)x_{\lambda }=A^{T}b,  \label{pt}
\end{equation}%
where $\lambda \in \mathbb{R}$ is a suitable parameter and $H$ is the
regularization matrix. The system (\ref{pt}) should have singular values
bounded away from $0$ in order to reduce the condition number and, at the
same time, its solution $x_{\lambda }$\ should be closed to the solution of
the original system.

For this kind of problem the method initially presented in this paper is
based on the shift and invert transformation%
\begin{equation}
Z=(A+\lambda I)^{-1},  \label{T1}
\end{equation}%
where $\lambda >0$ is a suitable parameter and $I$ is the identity matrix.
Provided that $\lambda $ is large enough, if $A\ $is positive definite ($%
F(A)\subset \mathbb{C}^{+}$, where $F(A)$ denotes the field of values) the
shift $A+\lambda I$, that represents the most elementary example of
regularization, has the immediate effect of moving the spectrum (that we
denote by $\sigma (A)$) away from $0$ so reducing the condition number.
Moreover, since%
\begin{equation*}
x=A^{-1}b=f(Z)b,
\end{equation*}%
where%
\begin{equation}
f(z)=\left( \frac{1}{z}-\lambda \right) ^{-1}=(1-\lambda z)^{-1}z,
\label{fz}
\end{equation}%
the idea is to solve the system $Ax=b$ by computing $f(Z)b$. For the
computation of $f(Z)b$, we use the standard Arnoldi method projecting the
matrix $Z$ onto the Krylov subspaces generated by $Z$ and $b$, that is $%
K_{m}(Z,b)=\mathrm{span}\{b,Zb,...,Z^{m-1}b\}$. By definition of $Z$ the
method is commonly referred to as the Restricted-Denominator (RD) rational
Arnoldi method \cite{Vanh}, \cite{Morno}.

Historically, a first attempt to reconstruct the solution from $x_{\lambda }$
that solves
\begin{equation}
\left( A+\lambda I\right) x_{\lambda }=b,  \label{st}
\end{equation}%
was proposed by Riley in \cite{Ri}. The algorithm is just based on the
approximation of $f(Z)$ by means of its Taylor series. Indeed we have%
\begin{equation}
A^{-1}b=\frac{1}{\lambda }\sum_{k=1}^{\infty }(\lambda Z)^{k}b,  \label{ser}
\end{equation}%
that leads to the recursion%
\begin{equation}
x_{k+1}=y+\lambda Zx_{k},\quad x_{0}=0,\quad y=Zb.  \label{ar}
\end{equation}%
It is easy to see that the method is equivalent to the \emph{iterative
improvement}%
\begin{eqnarray*}
\left( A+\lambda I\right) e_{k} &=&b-Ax_{k} \\
x_{k+1} &=&x_{k}+e_{k}
\end{eqnarray*}%
generally referred to as \emph{iterated Tikhonov regularization} or \emph{%
preconditioned Landweber iteration} (see e.g. \cite{Go}, \cite{HH}, \cite{KC}%
, \cite{KiCh}, \cite{Neu}). The main problem concerning this kind of
algorithms is that they can be extremely slow because the spectrum of $Z$
accumulates at $1/\lambda $ (cf. (\ref{T1}), (\ref{ser})). This, of course,
large values of $\lambda $, that is, when $A+\lambda I$ is well conditioned.
>From the point of view of the computation of function of matrices this is a
well known problem, i.e., the the computation by means of the Taylor series
generally provides poor results unless the spectrum of the matrix is close
to the expansion point. Indeed, from well known results of complex
approximation, the rate of convergence of a polynomial method for the
computation of a function of matrix depends on the position of the
singularity of the function, with respect to the location of the spectrum of
the matrix.

We also point out that, in \cite{BRSR}, the authors construct an improved
approximation via extrapolation with respect to the regularization
parameter, using the singular values representation of the solution.
Extrapolation techniques can also be applied to accelerate (\ref{ar}), as
suggested in \cite{BRED} and also indicated by Fasshauer in \cite{Fass}.

For problems in which the right hand side is affected by noise, instead of
working with the transformation (\ref{T1}) or implicitly with systems of
type (\ref{st}), we shall work with the standard regularization (\ref{pt})
and hence on the transformation%
\begin{equation*}
Z=(A^{T}A+\lambda L^{T}L)^{-1}.
\end{equation*}%
As we shall see, the subsequent Arnoldi-based algorithm for the
reconstruction of the exact solution will be almost identical to the one
based on (\ref{T1}), but the use of a regularization matrix $L$ different
from the identity allows to define methods less sensitive to perturbations
on the right hand side.

The paper is organized as follows. In Section \ref{due}, we describe the
Arnoldi method for the computation of $f(Z)b$ and, in Section \ref{tre}, we
present a theoretical a-priori error analysis. In Section \ref{quattro}, we
show an a-posteriori representation of the error. In Section \ref{cinque},
we analyze the choice of the parameter $\lambda $. Some numerical
experiments taken out from Hansen's Matlab toolbox on regularization \cite{H1,H2}%
, and from the theory of interpolation with radial basis functions are
presented in Section \ref{sei}. Finally, in Section \ref{sette}, we extend
our method to the Tikhonov regularization in its general form (\ref{pt})
showing also some tests with data affected by noise.

\section{The Arnoldi method for $f(Z)b$.}

\label{due}

For the construction of the subspaces $K_{m}(Z,b)$, the Arnoldi algorithm
generates an orthonormal sequence$\ \left\{ v_{j}\right\} _{j\geq 0}$, with $%
v_{1}=b/\left\Vert b\right\Vert $, such that $K_{m}(Z,b)=\mathrm{span}%
\left\{ v_{1},v_{2},...,v_{m}\right\} $ (here and below the norm used is
always the Euclidean norm). For every $m$ we have
\begin{equation}
ZV_{m}=V_{m}H_{m}+h_{m+1,m}v_{m+1}e_{m}^{T},  \label{cla}
\end{equation}%
where $V_{m}=\left[ v_{1},v_{2},...,v_{m}\right] $, $H_{m}$ is an upper
Hessenberg matrix with entries $h_{i,j}=v_{i}^{T}Zv_{j}$ and $e_{j}$ is the $%
j$-th vector of the canonical basis of \ $\mathbb{R}^{m}$. Formula (\ref{cla}%
) is just the matrix formulation of the algorithm.

The $m$-th Arnoldi approximation to $x=f(Z)b$ is defined as
\begin{equation*}
x_{m}=\left\Vert b\right\Vert V_{m}f(H_{m})e_{1}.
\end{equation*}%
Regarding the computation $f(H_{m})$, since the method is expected to
produce a good approximation of the solution in a relatively small number of
iterations, that is for $m\ll N$, one typically considers a certain rational
approximation to $f$, or the Schur-Parlett algorithm (see e.g. \cite[Chapter
11]{GV} or \cite{H}).

Denoting by $\Pi _{m-1}$ the vector space of polynomials of degree at most $%
m-1$, it can be seen that
\begin{equation}
x_{m}=\overline{p}_{m-1}(Z)b,  \label{pol}
\end{equation}%
where $\overline{p}_{m-1}\in $ $\Pi _{m-1}$ interpolates, in the Hermite
sense, the function $f$ at the eigenvalues of $H_{m}$ \cite{Saad2}.

As already mentioned, this kind of approach is commonly referred to as the
RD rational Arnoldi method since it is based on the use of single pole
rational forms of the type%
\begin{equation*}
R_{m-1}(x)=\frac{q_{m-1}(x)}{(x+a)^{m-1}},\mathbf{\quad }a\in \mathbb{R}%
,\quad q_{m-1}\in \Pi _{m-1},\quad m\geq 1,
\end{equation*}%
introduced and studied by N{\o }rsett in \cite{Nors} for the approximation
of the exponential function. In other words, with respect to $A$, formula (%
\ref{pol}) is actually a rational approximation.

It is worth noting that, at each step of the Arnoldi algorithm, we have to
compute the vectors $w_{j}=Zv_{j}$, $j\geq 1$, which leads to solve the
systems%
\begin{equation*}
(A+\lambda I)w_{j}=v_{j},\quad j\geq 1.
\end{equation*}%
Since $v_{1}=b/\left\Vert b\right\Vert $, the corresponding $w_{1}$ is just
the scaled solution of a regularized system (with the rough regularization $%
A\rightarrow A+\lambda I$). In this sense if $\lambda $ arises from the
standard techniques that seek for the optimal regularization parameter $%
\lambda _{opt}$ (L-curve, Generalized Cross Validation, etc.) this procedure
can be employed as a tool to improve the quality of the approximation $%
w_1\!\left\Vert b\right\Vert$. Anyway we shall see that, using the Arnoldi
algorithm, larger values for $\lambda $ are more reliable.

\section{Error analysis}

\label{tre}

The error $E_{m}:=x-x_{m}$ can be expressed and bounded in many ways (see
e.g. the recent paper \cite{BR} and the references therein). In any case,
however, the sharpness of the bound essentially depends on the amount of
information about the location of the field of values of $Z$, defined by%
\begin{equation*}
F(Z):=\left\{ \frac{x^{H}Zx}{x^{H}x},x\in \mathbb{C}^{N}\mathbf{\backslash }%
\left\{ 0\right\} \right\} .
\end{equation*}%
The bound we propose is based on the use of Faber polynomials. We need some
definitions and we refer to \cite{SL} or \cite{Wal} for a wide background of
what follows.

Let $\Omega $ be a compact and connected set of the complex plane. By the
Riemann mapping theorem there exists a conformal surjection
\begin{equation}
\psi :\overline{\mathbb{C}}\setminus \left\{ w:\left\vert w\right\vert \leq
1\right\} \rightarrow \overline{\mathbb{C}}\setminus \Omega ,\quad \psi
\left( \infty \right) =\infty ,\quad \psi ^{\prime }\left( \infty \right)
=\gamma ,  \label{3.1}
\end{equation}%
that has a Laurent expansion of the type%
\begin{equation*}
\psi (w)=\gamma w+c_{0}+\frac{c_{1}}{w}+\frac{c_{2}}{w^{2}}+\cdots
\end{equation*}%
The constant $\gamma $ is the capacity of $\Omega $. If $\Omega $ is an
ellipse or a line segment then $c_{i}=0$ for $i\geq 2$. Given a function $g$
analytic in $\Omega $, it is known that defining $p_{m-1}$ as the truncated
Faber series of exact degree $m-1$ with respect to $g$ and $\psi ,$ then $%
p_{m-1}$ provides an asymptotically optimal uniform approximation to $g$ in $%
\Omega $, that is%
\begin{equation}
\underset{m\rightarrow \infty }{\lim }\sup \left\Vert p_{m-1}-g\right\Vert
_{\Omega }^{1/m}=\underset{m\rightarrow \infty }{\lim }\sup \left\Vert
p_{m-1}^{\ast }-g\right\Vert _{\Omega }^{1/m},  \label{mc}
\end{equation}%
$\left\{ p_{m-1}^{\ast }\left( z\right) \right\} _{m\geq 1}$ being the
sequence of polynomials of best uniform approximation to $g$ in $\Omega $.
Property (\ref{mc}) is also called \emph{maximal convergence}. Let moreover $%
\phi :\overline{\mathbb{C}}\setminus \Omega \rightarrow \overline{\mathbb{C}}%
\setminus \left\{ w:\left\vert w\right\vert \leq 1\right\} $ be the inverse
of $\psi $. For any $r>1,$ let $\Gamma _{r}$ be the equipotential curve
\begin{equation*}
\Gamma _{r}:=\left\{ z:\left\vert \phi \left( z\right) \right\vert
=r\right\} ,
\end{equation*}%
and let us denote by $\Omega _{r}$ the bounded domain with boundary $\Gamma
_{r}$. Let $\widehat{r}>1$ be the largest number such that $g$ is analytic
in $\Omega _{r}$ for each $\gamma <r<\widehat{r}$ and has a singularity on $%
\Gamma _{\widehat{r}}$. Then, it is known that the rate of convergence of
the sequence $\left\{ p_{m-1}\left( z\right) \right\} _{m\geq 1}$ is given by%
\begin{equation}
\underset{m\rightarrow \infty }{\lim }\sup \left\Vert p_{m-1}-g\right\Vert
_{\Omega }^{1/m}=\frac{1}{\widehat{r}}.  \label{mp}
\end{equation}%
For this reason we know that superlinear convergence is only attainable for
entire functions, where asymptotically one can set $\widehat{r}:=m$. In
order to derive error bounds for the computation of $f(Z)b$ we need the
following classical result

\begin{theorem}
\textrm{\cite{Ellac}} Let $\Omega $ be a compact and convex subset such that
$g$ is analytic in $\Omega $. For $1<r<\widehat{r}$ the following bound
holds
\begin{equation}
\left\Vert p_{m-1}-g\right\Vert _{\Omega }\leq 2\left\Vert g\right\Vert
_{\Gamma _{r}}\frac{\displaystyle \left( \frac{1}{r}\right) ^{m}}{%
\displaystyle 1-\frac{1}{r}}.  \label{4.3b}
\end{equation}
\end{theorem}

Using the above theorem, for our function $f(z)=z/(1-\lambda z)$, singular
at $1/\lambda $, we can state the

\begin{proposition}
\label{p1}Assume that $\Omega $ is an ellipse of the complex plane,
symmetric with respect to the real axis with associated conformal mapping $%
\psi (w)=\gamma w+c_{0}+c_{1}/w$. Assume that $\psi (1)<1/\lambda $ and let $%
\widehat{r}$ be such that $\psi (\widehat{r})=1/\lambda $. Let moreover $%
\overline{m}$ be the smallest integer such that%
\begin{equation*}
\frac{\widehat{r}}{\overline{m}+1}<\widehat{r}-1.
\end{equation*}%
Then for $m\geq \overline{m}$%
\begin{equation}
\left\Vert p_{m-1}-f\right\Vert _{\Omega }\leq \frac{2\, e \,\overline{m}\,
\widehat{r}}{\overline{m}(\widehat{r}-1)-1}\frac{1}{\lambda ^{2}\psi
^{\prime }(\widehat{r})}\frac{m+1}{\widehat{r}^{m}},  \label{f1e}
\end{equation}%
and for $m<\overline{m}$%
\begin{equation}
\left\Vert p_{m-1}-f\right\Vert _{\Omega }\leq \frac{4}{\lambda ^{2}\left(
\widehat{r}-1\right) \psi ^{\prime }(\widehat{r})}\left( \frac{2}{\widehat{r}%
+1}\right) ^{m}\frac{\widehat{r}+1}{\widehat{r}-1}.  \label{f2e}
\end{equation}
\end{proposition}

\begin{proof}
Let $r=\widehat{r}-\varepsilon $, with $0<\varepsilon <\widehat{r}-1$. By
the properties of $\Omega $, we have%
\begin{equation*}
\left\Vert f\right\Vert _{\Gamma _{r}}=\frac{\psi (r)}{1-\lambda \psi (r)},
\end{equation*}%
and, by direct computation%
\begin{equation*}
\psi (r)=\psi (\widehat{r})-\gamma \varepsilon +\frac{c_{1}\varepsilon }{(%
\widehat{r}-\varepsilon )\widehat{r}}.
\end{equation*}%
Hence using $\psi (\widehat{r})=1/\lambda $ we find%
\begin{eqnarray*}
\left\Vert f\right\Vert _{\Gamma _{r}} &\leq &\frac{\psi (\widehat{r})}{%
1-\lambda \left( \psi (\widehat{r})-\gamma \varepsilon +\frac{\displaystyle %
c_{1}\varepsilon }{\displaystyle (\widehat{r}-\varepsilon )\widehat{r}}%
\right) }, \\
&=&\frac{1}{\lambda ^{2}\varepsilon \left( \gamma -\frac{\displaystyle c_{1}%
}{\displaystyle (\widehat{r}-\varepsilon )\widehat{r}}\right) }, \\
&\leq &\frac{1}{\lambda ^{2}\varepsilon \psi ^{\prime }(\widehat{r})}.
\end{eqnarray*}%
By (\ref{4.3b}), we thus obtain%
\begin{equation}
\left\Vert p_{m-1}-f\right\Vert _{\Omega }\leq \frac{2}{\lambda
^{2}\varepsilon \psi ^{\prime }(\widehat{r})}\frac{1}{\left( \widehat{r}%
-\varepsilon \right) ^{m}}\frac{1}{\displaystyle 1-\frac{1}{\widehat{r}%
-\varepsilon }}.  \label{ers}
\end{equation}%
Now setting
\begin{equation}
\varepsilon =\frac{\widehat{r}}{m+1},  \label{ep}
\end{equation}%
since this value minimizes%
\begin{equation*}
\frac{1}{\varepsilon \left( \widehat{r}-\varepsilon \right) ^{m}},
\end{equation*}%
let $\overline{m}$ be the smallest positive integer such that%
\begin{equation*}
\frac{\widehat{r}}{\overline{m}+1}<\widehat{r}-1.
\end{equation*}%
By inserting (\ref{ep}) into (\ref{ers}) and using%
\begin{equation*}
\frac{1}{\displaystyle 1-\frac{1}{\widehat{r}-\varepsilon }}\leq \frac{%
\overline{m}\widehat{r}}{\overline{m}(\widehat{r}-1)-1},
\end{equation*}%
we find (\ref{f1e}). For $m<\overline{m}$ we can take for instance%
\begin{equation}
\varepsilon =\frac{\widehat{r}-1}{2}.  \label{ep2}
\end{equation}%
Substituting (\ref{ep2}) into (\ref{ers}) we obtain (\ref{f2e}).
\end{proof}

\begin{remark}
Note that the assumption $\psi (1)<1/\lambda $ in Proposition \ref{p1} just
means that the ellipse is strictly on the left of the singularity of $f$.
\end{remark}

Regarding the field of values of $Z$, $F(Z)$, it is well known that it is
convex, that $\sigma (Z)\subset F(Z)$, and that $F(H_{m})\subseteq F(Z)$
(where $H_m$ is defined in Section \ref{due}). Of course if $F(A)\subset
\mathbb{C}^{+}$ ($A$ is positive definite) then $F(Z)\subset \{z\in C:0<%
\mathrm{Re}(z)<1/\lambda \}$ and the corresponding $f$ is analytic in $F(Z)$%
. Using these properties we can state the following result

\begin{theorem}
\label{t1}Assume that $F(A)\subset \mathbb{C}^{+}$. Let $\Omega $ be an
ellipse (with associated conformal mapping $\psi $, and inverse $\phi $)
symmetric with respect to the real axis and such that $F(Z)\subseteq \Omega $
with $f$ analytic in $\Omega $. Then, for $m$ large enough, we have%
\begin{equation*}
\left\Vert E_{m}\right\Vert \leq 4\,e\,C\frac{\widehat{r}}{\widehat{r}-1}%
\frac{1}{\psi ^{\prime }(\widehat{r})}\,K\,\frac{m+1}{\widehat{r}^{m}},
\end{equation*}%
where $K=1/\lambda ^{2}$, $\widehat{r}=\phi (1/\lambda )$, and $C=$ $11.08$ (%
$C=1$ if $A$ is symmetric).
\end{theorem}

\begin{proof}
Using the properties of the Arnoldi algorithm, we know that for every $%
p_{m-1}\in \Pi _{m-1}$,
\begin{equation}
V_{m}p_{m-1}(H_{m})e_{1}=p_{m-1}(Z)b.  \label{re1}
\end{equation}%
Hence, from (\ref{re1}), it follows that, for $m\geq 1$ and for every $%
p_{m-1}\in \Pi _{m-1}$,
\begin{equation}
E_{m}=x-x_{m}=f(Z)b-p_{m-1}(Z)b-V_{m}(f(H_{m})-p_{m-1}(H_{m}))e_{1}.
\label{e1}
\end{equation}%
Since $\left\Vert V_{m}\right\Vert =1$ we have (see \cite{Cru})
\begin{equation}
\left\Vert E_{m}\right\Vert \leq 2C\left\Vert p_{m-1}-f\right\Vert _{F(Z)}.
\label{errf}
\end{equation}%
Therefore taking $p_{m-1}$ as the $\left( m-1\right) $-th truncated Faber
(Chebyshev) series, the result follows from Proposition \ref{p1} since $%
F(Z)\subseteq \Omega $.
\end{proof}

\begin{remark}
By (\ref{e1}), if both $Z$ and $H_{m}$ are diagonalizable then $C$ in (\ref%
{errf}) is a constant depending on the condition number of the
diagonalization matrices and $\Omega $ can be taken as an ellipse containing
$\sigma (A)$.
\end{remark}

Theorem \ref{t1} is surely important from a theoretical point of view since
it states that the Arnoldi algorithm produces asymptotically optimal
approximations. However, if we consider for simplicity the symmetric case,
we can also understand that it cannot be used to suggest the choice of $%
\lambda$.

Indeed, let $\lambda _{1}\gtrsim 0$ and $\lambda _{N}$ be respectively the
smallest and the largest eigenvalues $A$. Then $F(A)=[\lambda _{1},\lambda
_{N}]$ and
\begin{equation*}
\displaystyle F(Z)=\left[ \frac{1}{\lambda _{N}+\lambda },\frac{1}{\lambda
_{1}+\lambda }\right] =:I_{\lambda }.
\end{equation*}%
In this case, by (\ref{errf}) we have%
\begin{equation*}
\left\Vert E_{m}\right\Vert \leq 2\max_{I_{\lambda }}\left\vert
f(z)-p_{m-1}(z)\right\vert .
\end{equation*}%
As already mentioned, the conformal mapping $\psi $ associated to $%
I_{\lambda }$ takes the form%
\begin{equation}
\psi (w)=\gamma w+c_{0}+\frac{c_{1}}{w}  \label{conf}
\end{equation}%
where%
\begin{eqnarray}
\gamma  &=&\frac{1}{4}\left( \frac{1}{\lambda _{1}+\lambda }-\frac{1}{%
\lambda _{N}+\lambda }\right) =\frac{1}{4}\frac{\lambda _{N}-\lambda _{1}}{%
\left( \lambda _{1}+\lambda \right) (\lambda _{N}+\lambda )},  \notag \\
c_{0} &=&\frac{1}{2}\left( \frac{1}{\lambda _{1}+\lambda }+\frac{1}{\lambda
_{N}+\lambda }\right) =\frac{1}{2}\frac{\lambda _{N}+\lambda _{1}+2\lambda }{%
\left( \lambda _{1}+\lambda \right) (\lambda _{N}+\lambda )},  \label{c00} \\
c_{1} &=&\gamma \text{.}  \notag
\end{eqnarray}%
For $r>1$, $\Omega _{r}$ is the confocal ellipse (foci in $\displaystyle%
\frac{1}{\lambda _{N}+\lambda }$ and $\displaystyle\frac{1}{\lambda
_{1}+\lambda }$) described by $\psi (re^{i\theta })$, $0\leq \theta <2\pi $.
Since $f(z)$ is singular at $1/\lambda $, $\widehat{r}$ is the solution ($>1$%
) of
\begin{equation}
\gamma \widehat{r}+c_{0}+\frac{\gamma }{\widehat{r}}=\frac{1}{\lambda }
\label{rs}
\end{equation}%
that is%
\begin{equation}
\widehat{r}=u+\sqrt{u^{2}-1},  \label{r1}
\end{equation}%
where%
\begin{equation}
u=\frac{2\lambda _{1}\lambda _{N}}{\lambda (\lambda _{N}-\lambda _{1})}+%
\frac{\lambda _{N}+\lambda _{1}}{\lambda _{N}-\lambda _{1}}.  \label{r2}
\end{equation}%
Thus, $\widehat{r}$ monotonically decreases with respect to $\lambda $ and $%
\widehat{r}\rightarrow \infty $ for $\lambda \rightarrow 0$.

The above arguments simply show that the error analysis does not take into
account of the computational problems in the inversion of $A+\lambda I$ for $%
\lambda \approx 0$. The method is very fast for $\lambda \approx 0$ because,
at each step, we are inverting something very close to the original operator
$A$. In order to derive a more useful estimate one should modify the above
analysis imposing in some way the requirement $\lambda \gg \lambda _{1}$. In
some sense this will be done in Section \ref{cinque} where we consider the
conditioning in the computation of $f(Z)b$ that is obviously closely related
to the rate of convergence of any iterative method.

\section{A-posteriori error representation}

\label{quattro}

By a result on Pad\'{e}--type approximation proved in \cite{CBout}, we know
that the Hermite interpolation polynomial of the function
\begin{equation*}
g(s)=\frac{1}{1-st}
\end{equation*}%
at the zeros of any polynomial $\nu _{m}$ of exact degree $m$ in $s$ is
given by%
\begin{equation*}
R_{m-1}(s)=\frac{1}{1-st}\left( 1-\frac{\nu _{m}(s)}{\nu _{m}(t^{-1})}%
\right) .
\end{equation*}%
Setting $\lambda =t^{-1}$, we have that%
\begin{equation*}
f(\xi )=\frac{1}{\xi ^{-1}-\lambda }=-\lambda ^{-1}g\left( \xi ^{-1}\right) ,
\end{equation*}%
and so%
\begin{equation}
-\lambda ^{-1}R_{m-1}(\xi ^{-1})=\frac{1}{1-\xi ^{-1}\lambda ^{-1}}\left( 1-%
\frac{\nu _{m}(\xi ^{-1})}{\nu _{m}(\lambda )}\right)  \label{s}
\end{equation}%
interpolates $f(\xi )$. By (\ref{pol}) let $\overline{p}_{m-1}\in $ $\Pi
_{m-1}$ be the polynomial that interpolates, in the Hermite sense, the
function $f(z)$ at the eigenvalues of $H_{m}$, $\xi _{1},...,\xi _{m^{\prime
}}$, $m^{\prime }\leq m$, with multiplicity $k_{i}$, $i=1,...,m^{\prime }$%
.Then
\begin{equation*}
\overline{p}_{m-1}^{(j)}(\xi _{i})=-\lambda ^{-1}R_{m-1}^{(j)}(\xi
_{i}^{-1})=f^{(j)}(\xi _{i}),\quad 1\leq i\leq m^{\prime },\ 0\leq j\leq
k_{i}-1.
\end{equation*}%
By (\ref{s}) and using the above relation is it easy to see that $\nu
_{m}(s)=\det (sI-H_{m}^{-1})$. In this way, by direct computation,%
\begin{eqnarray}
x_{m} &=&\overline{p}_{m-1}(Z)b,  \notag \\
&=&A^{-1}b-A^{-1}\left( \frac{\nu _{m}(Z^{-1})}{\nu _{m}(\lambda )}\right) b.
\label{e2}
\end{eqnarray}%
Since, of course, $A^{-1}$ and $Z^{-1}$ commute, we find%
\begin{equation*}
\frac{\left\Vert x_{m}-x\right\Vert }{\left\Vert x\right\Vert }\leq \frac{%
\left\Vert \nu _{m}(A+\lambda I)\right\Vert }{\left\vert \nu _{m}(\lambda
)\right\vert }.
\end{equation*}

A posteriori error estimate can be derived in this way. Since%
\begin{eqnarray*}
\nu _{m}(s) &=&\det (sI-H_{m}^{-1}), \\
&=&\frac{s^{m}\det (H_{m}-s^{-1}I)}{\det H_{m}},
\end{eqnarray*}%
defining $q_{m}(\xi )=\det (H_{m}-\xi I)$, we have
\begin{equation}
\frac{\left\Vert x_{m}-x\right\Vert }{\left\Vert x\right\Vert }\leq \frac{%
\left\Vert (A+\lambda I)^{m}q_{m}(Z)\right\Vert }{\lambda ^{m}\left\vert
q_{m}(\lambda ^{-1})\right\vert }.  \label{ape}
\end{equation}

It is worth noting that, using the relation
\begin{equation*}
q_{m}(Z)b=\left( \prod\nolimits_{j=1}^{m}h_{j+1,j}\right) v_{m+1},
\end{equation*}%
(see \cite{Morno}), we obtain from (\ref{e2})%
\begin{equation*}
\left\Vert x_{m}-x\right\Vert =\frac{\left(
\prod\nolimits_{j=1}^{m}h_{j+1,j}\right) }{\lambda ^{m}\left\vert
q_{m}(\lambda ^{-1})\right\vert }\left\Vert A^{-1}(A+\lambda
I)^{m}v_{m+1}\right\Vert ,
\end{equation*}%
which proves the convergence in a finite number $m^{\ast }\leq N$ of steps
of the method in exact arithmetics. Note that by (\ref{e2}) the
corresponding $\nu _{m^{\ast }}$ is the minimal polynomial of $A+\lambda I$
for the vector $b$.

\section{The choice of $\protect\lambda $}

\label{cinque}

\begin{figure}[h]
\begin{center}
\includegraphics[width=15cm]{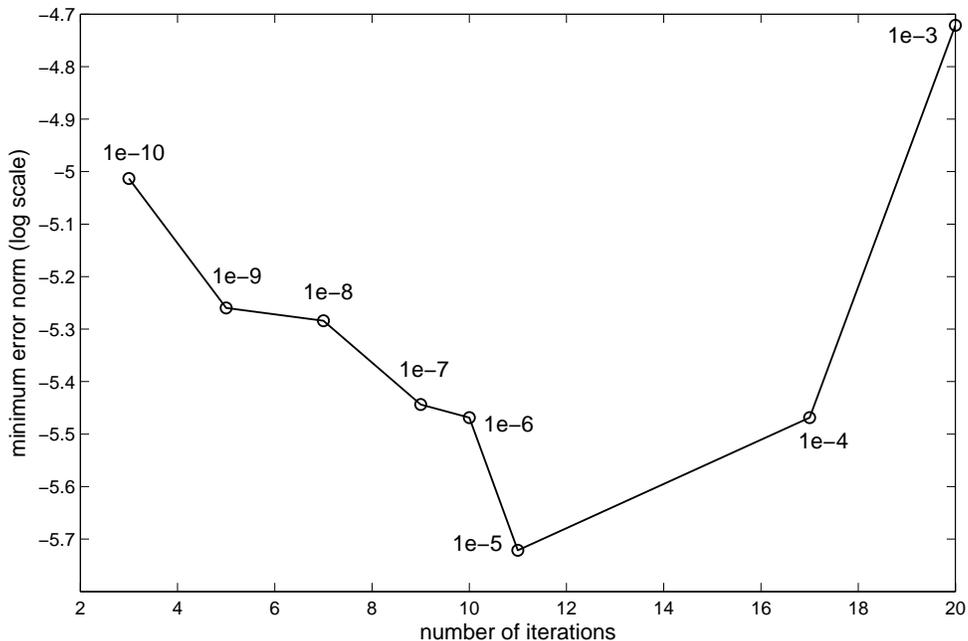}
\end{center}
\caption{BAART(40) - Minimum attained error with respect to the number of
iterations for different values of $\protect\lambda$. }
\label{BBART40}
\end{figure}

As already mentioned, the arguments of Section \ref{tre} reveal that the
standalone error analysis of the computation of $f(Z)b$ is not reliable to
suggest the choice of $\lambda $, since $\kappa (Z)\rightarrow \kappa (A)$
as $\lambda \rightarrow 0$ ($\kappa (\cdot )$ denoting the standard
condition number of a matrix). In other words, it does not take into account
that, at each step, we need to solve a system with the matrix $A+\lambda I$.
At the same time, focusing the attention on the accuracy (so neglecting the
rate of convergence) one could expect that "large" values of $\lambda $
should allow an improvement of it, since the linear systems with $A+\lambda
I $ would be solved more accurately. The numerical experiments show that
this is not true, as shown in Fig. \ref{BBART40}, where we consider the
problem BAART, taken out from the Hansen's Matlab toolbox \texttt{Regtools} (see
\cite{H1} and \cite{H2}).

Indeed the diagram of Fig. \ref{BBART40} represents the standard situation,
that is, increasing $\lambda $, we have a loss of accuracy. The behavior on
the leftmost part of the diagram is clear since it is due to the
conditioning of $Z$ for $\lambda $ small. On the rightmost part we have
again a loss of accuracy but now it depends on the numerical instability in
the computation of $f(Z)$ for $\lambda $ large (the problem can be easily
observed even working scalarly). This observation leads us to consider the
conditioning in the computation of $f(Z)b$ for having a good strategy to
define $\lambda $.

The absolute and the relative condition number for the computation of $g(X)$
where $g$ is a given function and $X$ a square matrix are given by (cf. \cite%
{H} Chapter 3)%
\begin{eqnarray}
\kappa _{a}(g,X) &=&\lim_{\varepsilon \rightarrow 0}\sup_{\left\Vert
E\right\Vert \leq \varepsilon }\frac{\left\Vert g(X+E)-g(X)\right\Vert }{%
\varepsilon },  \label{ka} \\
\kappa _{r}(g,X) &=&\kappa _{a}(g,X)\frac{\left\Vert X\right\Vert }{%
\left\Vert g(X)\right\Vert },  \label{kr}
\end{eqnarray}%
and these definitions imply that%
\begin{equation*}
\left\Vert g(X+E)-g(X)\right\Vert \leq \kappa _{a}(g,X)\left\Vert
E\right\Vert +O(\left\Vert E\right\Vert ^{2}).
\end{equation*}

\begin{proposition}
For the function $f(z)=(1-\lambda z)^{-1}z$ we have the bound%
\begin{equation}
\kappa _{r}(f,Z)\leq \frac{\left\Vert (I-\lambda Z)^{-2}\right\Vert
\left\Vert Z\right\Vert }{\left\Vert (Z^{-1}-\lambda I)^{-1}\right\Vert }.
\label{bc}
\end{equation}
\end{proposition}

\begin{proof}
In order to derive first the absolute condition number we have%
\begin{eqnarray*}
f(Z+E)-f(Z) &=&\left[ (Z+E)^{-1}-\lambda I\right] ^{-1}-(Z^{-1}-\lambda
I)^{-1}, \\
&=&\left[ (I+Z^{-1}E)^{-1}Z^{-1}-\lambda I\right] ^{-1}-(Z^{-1}-\lambda
I)^{-1}, \\
&=&\left[ Z^{-1}-\lambda I+\Lambda (Z,E)\right] ^{-1}-(Z^{-1}-\lambda
I)^{-1},
\end{eqnarray*}%
where%
\begin{equation*}
\Lambda (Z,E):=\sum_{k=1}^{\infty }(-1)^{k}(Z^{-1}E)^{k}Z^{-1}.
\end{equation*}%
Hence%
\begin{eqnarray}
f(Z+E)-f(Z) &=&\left[ I+(Z^{-1}-\lambda I)^{-1}\Lambda (Z,E)\right]
^{-1}(Z^{-1}-\lambda I)^{-1}-(Z^{-1}-\lambda I)^{-1},  \notag \\
&=&\sum\nolimits_{j=0}^{\infty }(-1)^{j}(Z^{-1}-\lambda I)^{-j}\Lambda
(Z,E)^{j}(Z^{-1}-\lambda I)^{-1}-(Z^{-1}-\lambda I)^{-1},  \label{fd}
\end{eqnarray}%
and finally%
\begin{equation*}
\left\Vert f(Z+E)-f(Z)\right\Vert \leq \left\Vert (Z^{-1}-\lambda
I)^{-1}Z^{-1}EZ^{-1}(Z^{-1}-\lambda I)^{-1}\right\Vert +O(\left\Vert
E\right\Vert ^{2}),
\end{equation*}%
so that%
\begin{equation*}
\kappa _{a}(f,Z)\leq \left\Vert (I-\lambda Z)^{-2}\right\Vert ,
\end{equation*}%
that proves (\ref{bc}) using (\ref{kr}) and the definition of $f(z)$. Note
that by (\ref{fd})%
\begin{equation*}
L(Z,E):=(I-\lambda Z)^{-1}E(I-\lambda Z)^{-1}
\end{equation*}%
is the Fr\'{e}chet derivative of $f$ at $Z$ applied to $E$.
\end{proof}

This Proposition simply shows that the problem is well conditioned for $%
\lambda \rightarrow 0$ and ill conditioned for $\lambda \gg 0$, that matches
with the error analysis of Section \ref{tre}. Of course the situation is
opposite to what happens for the solution of the linear systems with $%
A+\lambda I$ during the Arnoldi process. Therefore the idea, confirmed by
many numerical experiments, is to define $\lambda $ such that $\kappa
_{r}(f,Z)\approx \kappa (A+\lambda I)$, that is, to consider the bound (\ref%
{bc}) and solve the equation%
\begin{equation*}
\frac{\left\Vert (I-\lambda Z)^{-2}\right\Vert \left\Vert Z\right\Vert }{%
\left\Vert (Z^{-1}-\lambda I)^{-1}\right\Vert }=\left\Vert (A+\lambda
I)\right\Vert \left\Vert (A+\lambda I)^{-1}\right\Vert .
\end{equation*}

In the SPD case everything becomes clear since we have%
\begin{eqnarray*}
\frac{\left\Vert (I-\lambda Z)^{-2}\right\Vert \left\Vert Z\right\Vert }{%
\left\Vert (Z^{-1}-\lambda I)^{-1}\right\Vert } &=&\frac{\lambda +\lambda
_{1}}{\lambda _{1}} \\
\left\Vert (A+\lambda I)\right\Vert \left\Vert (A+\lambda I)^{-1}\right\Vert
&=&\frac{\lambda _{N}+\lambda }{\lambda _{1}+\lambda }
\end{eqnarray*}%
that for $\lambda _{1}\rightarrow 0$ leads to%
\begin{equation*}
\lambda =\sqrt{\lambda _{1}\lambda _{N}}+O(\lambda _{1}).
\end{equation*}

\begin{remark}
If the underlying operator is bounded then one may consider the approximation%
\begin{equation*}
\sqrt{\lambda _{1}\lambda _{N}}\approx \frac{1}{\sqrt{\kappa (A)}}\quad
\text{for }\lambda _{1}\rightarrow 0.
\end{equation*}
\end{remark}

\begin{remark}
In the SPD case, taking $\lambda ^{\ast }=\sqrt{\lambda _{1}\lambda _{N}}$
and putting it into (\ref{r1})-(\ref{r2}), we find that the asymptotic
convergence factor of the method is given by%
\begin{equation*}
\left\Vert E_{m}\right\Vert ^{1/m}\rightarrow \frac{1}{\widehat{r}}=\frac{%
\lambda _{N}^{1/4}-\lambda _{1}^{1/4}}{\lambda _{N}^{1/4}+\lambda _{1}^{1/4}}%
=\frac{\kappa (A)^{1/4}-1}{\kappa (A)^{1/4}+1}.
\end{equation*}
\end{remark}

\begin{remark}
The choice of $\lambda ^{\ast }$ has another interesting meaning. Indeed,
let us consider the problem of the computation of $g(A)b$ with $g$ singular
only at 0 and $A$ SPD. Using the transformation $z=\left( a+\lambda \right)
^{-1}$ (cf. (\ref{T1})), if the corresponding $g^{\ast }(z)=g(z^{-1}-\lambda
)$ has a non-removable singularity at 0, then the optimal choice of $\lambda
$ is given by solving the equation%
\begin{equation}
c_{0}=\frac{1}{2\lambda }  \label{c0}
\end{equation}%
(cf. (\ref{conf}) and (\ref{c00})), that is, the midpoint of $[0,1/\lambda ]$
must be equal to the midpoint of $I_{\lambda }$, because in this way we have
simultaneously $\psi (-\widehat{r})=0$ and $\psi (\widehat{r})=1/\lambda $.
A straightforward computation shows that solving (\ref{c0}) leads exactly to
$\lambda ^{\ast }$. For instance, in \cite{Mo} the author uses the RD
Arnoldi method to compute $\sqrt{A}b$ and obtains the same result even if
following a different approach.
\end{remark}

\begin{remark}
The condition number of $A+\lambda ^{\ast }I$ is given by%
\begin{equation*}
\kappa (A+\lambda ^{\ast }I)=\frac{\lambda _{N}+\sqrt{\lambda _{1}\lambda
_{N}}}{\lambda _{1}+\sqrt{\lambda _{1}\lambda _{N}}}=\sqrt{\frac{\lambda _{N}%
}{\lambda _{1}}}=\sqrt{\kappa (A)}.
\end{equation*}
\end{remark}

In the nonsymmetric case, the analysis is a bit more difficult but many
numerical experiments have shown that just having information on the
conditioning of $A$, the choice $\lambda \approx \kappa (A)^{-1/2}$ is
generally satisfactory, that is, we are rather close to the minimum of a
curve similar to the one of Fig. \ref{BBART40}. For very ill-conditioned
problems we suggest to define $\lambda $ a bit larger, say in the range $%
10\kappa (A)^{-1/2}\div 100\kappa (A)^{-1/2}$, since the errors generated by
the solution of the linear systems might be much larger than the machine
precision.

\section{Numerical experiments}

\label{sei}

In order to test the efficiency of our method, that from now on we denote by
RA (Rational Arnoldi), we consider here some numerical experiments where we
compare it with other classical iterative solvers. The RA method have have
been implemented in Matlab following the line of Algorithm \ref{alg_CMP}
described below.

\begin{algorithm}[ht]
\begin{algo}
\REQUIRE $A\in {\mathbb{R}}^{N\times N}\!\!,\;b\in {\mathbb{R}}^{N}, \lambda \in {\mathbb{R}}$
\DEFINE $f=(1-\lambda z)^{-1}z$
\STATE[3:] \If $ (A+\lambda I)$ is \textsc{SPD},
\Then\textbf{Compute} $L$ s.t. $(A+\lambda
I)=L\,L^{T}$
\STATE[~] \Else {\bf Compute} $L,U$ s.t. $(A+\lambda I)=L\,U$, %
\Endif\STATE[4:] $v_1\leftarrow b/\Vert b\Vert, V_{1}\leftarrow [v_1  ]$
\FOR[5:] $\!\!m=1,2,\ldots $
\Do
\UPDATE[5.1:] $H_{m}\in {\mathbb{R}}^{m\times m}$ by
Arnoldi's algorithm
\STATE[~] {\bf Remark:} In the Arnoldi's algorithm, we
compute $w_{m}=Zv_{m}$
\STATE[~] solving $(A+\lambda I)w_{m}=v_{m}$, that is
$w_{m}=U^{-1}L^{-1}v_{m}$ or $w_{m}=(L^{T})^{-1}L^{-1}v_{m}$.
\COMPUTE[5.2:]
$f(H_{m})$ by Schur-Parlett algorithm
\STATE[5.3:] $x_{m}\leftarrow \Vert
b\Vert V_{m}f(H_{m})\;e_{1}$
\OUTPUT[5.4:] $x_{m}$, approximation of $%
f(Z)b=A^{-1}b$
\UPDATE[5.5:] $V_{m+1}=[v_{1},\ldots ,v_{m+1}]\in {\mathbb{R}}%
^{N\times (m+1)}$ orthonormal basis for
\STATE[~] $K_{m+1}(Z,b)$, by
Arnoldi's algorithm
\ENDFOR[~]
\end{algo}

\caption{-  RA Algorithm for solving $Ax=b$.} \label{alg_CMP}
\end{algorithm}

It is worth noting that we make use of\ the LU (or Cholesky) factorization
to solve the linear system at each step. The reason is to reduce the
computational cost since the factorization is computed only once at the
beginning, taking also into account that $A+\lambda I$ should be relatively
well conditioned. Anyway, for large scale non-sparse problems an iterative
approach producing an inner-outer iteration should be considered.

We consider four classical test problems taken out from Hansen's Matlab toolbox
\texttt{Regtools}, GRAVITY, FOXGOOD, SHAW and BAART. These discrete linear
problems arise from the discretization of Fredholm integral equations of the
first kind. In all experiments, we consider a noise-free right hand side,
that is, we define $b=Ax$. The numerical results have been obtained with
Matlab 7.9, on a single processor computer Intel Core2 Duo T5800.

Tables \ref{TABLE1} and \ref{TABLE2} below summarize the results. For
comparison, we consider the codes ART, CGLS, LSQR\_B and MR2 taken out from
Hansen's toolbox, CG, GMRES and MINRES that are resident Matlab functions,
and Riley's method. The number between parentheses beside the name of the
test is the dimension of the system. In all tests $\lambda _{RA}$ and $%
\lambda _{Riley}$ denote the chosen values of the parameters for the RA and
Riley's method respectively. Since no general indication about the choice of
the parameter for Riley's method is available in the literature, in all
experiments we heuristically select a nearly best one. In the tables we
consider the minimum attained error norm \textit{err}, the corresponding
residual \textit{res} and the number of iterations \textit{nit}. Each method
was stopped when the number of iterations reaches the dimension of the
system. The missing numbers are due to the structure of the coefficient
matrix (symmetric, SPD, and so on).

\begin{table}[th]
\begin{center}
\begin{tabular}{l|c|c|l|c|c|l|}
& \multicolumn{3}{|c|}{GRAVITY(100)} & \multicolumn{3}{|c|}{FOXGOOD(80)} \\
\hline\hline
$\lambda _{RA}$, $\lambda _{\mathit{Riley}}$ & \multicolumn{3}{|c}{\textbf{%
1e-9}, 1e-11} & \multicolumn{3}{|c|}{\textbf{1e-8}, 1e-10} \\ \hline
& \textit{err} & \textit{res} & \textit{nit} & \textit{err} & \textit{res} &
\textit{nit} \\ \hline
\textbf{RA} & \multicolumn{1}{|l|}{\textbf{1.6e-5}} & \multicolumn{1}{|l|}{%
\textbf{8.1e-9}} & \textbf{2} & \multicolumn{1}{|l|}{\textbf{6.8e-7}} &
\multicolumn{1}{|l|}{\textbf{2.9e-10}} & \textbf{5} \\
CG & \multicolumn{1}{|l|}{1.7e-4} & \multicolumn{1}{|l|}{7.5e-11} & 96 &
\multicolumn{1}{|l|}{} & \multicolumn{1}{|l|}{} &  \\
ART & \multicolumn{1}{|l|}{8.4e-2} & \multicolumn{1}{|l|}{5.8e-3} & 100 &
\multicolumn{1}{|l|}{2.3e-3} & \multicolumn{1}{|l|}{8.8e-6} & 80 \\
CGLS & \multicolumn{1}{|l|}{} & \multicolumn{1}{|l|}{} &  &
\multicolumn{1}{|l|}{6.3e-6} & \multicolumn{1}{|l|}{9.6e-14} & 80 \\
LSQR\_B & \multicolumn{1}{|l|}{1.7e-3} & \multicolumn{1}{|l|}{2.0e-8} & 100
& \multicolumn{1}{|l|}{2.9e-6} & \multicolumn{1}{|l|}{1.1e-14} & 80 \\
MR2 & \multicolumn{1}{|l|}{1.9e-3} & \multicolumn{1}{|l|}{2.3e-8} & 66 &
\multicolumn{1}{|l|}{2.3e-6} & \multicolumn{1}{|l|}{1.6e-15} & 57 \\
MINRES & \multicolumn{1}{|l|}{1.8e-4} & \multicolumn{1}{|l|}{4.6e-11} & 100
& \multicolumn{1}{|l|}{2.0e-5} & \multicolumn{1}{|l|}{1.6e-15} & 80 \\
RILEY & \multicolumn{1}{|l|}{1.3e-3} & \multicolumn{1}{|l|}{8.0e-11} & 2 &
\multicolumn{1}{|l|}{6.3e-6} & \multicolumn{1}{|l|}{5.2e-10} & 2%
\end{tabular}%
\end{center}
\caption{Results for GRAVITY and FOXGOOD. }
\label{TABLE1}
\end{table}

\begin{table}[th]
\begin{center}
\begin{tabular}{l|c|c|l|c|c|l|}
& \multicolumn{3}{|c|}{SHAW(64)} & \multicolumn{3}{|c|}{BAART(120)} \\
\hline\hline
$\lambda _{RA}$, $\lambda _{\mathit{Riley}}$ & \multicolumn{3}{|c}{\textbf{%
1e-9}, 1e-10} & \multicolumn{3}{|c|}{\textbf{1e-8}, 1e-10} \\ \hline
& \textit{err} & \textit{res} & \multicolumn{1}{|c|}{\textit{nit}} & \textit{%
err} & \textit{res} & \textit{nit} \\ \hline
\textbf{RA} & \multicolumn{1}{|l|}{\textbf{3.3e-3}} & \multicolumn{1}{|l|}{%
\textbf{2.0e-7}} & \textbf{7} & \multicolumn{1}{|l|}{\textbf{8.3e-6}} &
\multicolumn{1}{|l|}{\textbf{1.3e-8}} & \textbf{6} \\
GMRES & \multicolumn{1}{|l|}{} & \multicolumn{1}{|l|}{} &  &
\multicolumn{1}{|l|}{9.6e-6} & \multicolumn{1}{|l|}{1.4e-15} & 15 \\
ART & \multicolumn{1}{|l|}{7.7e-1} & \multicolumn{1}{|l|}{6.8e-2} & 64 &
\multicolumn{1}{|l|}{3.4e-1} & \multicolumn{1}{|l|}{2.7e-2} & 120 \\
CGLS & \multicolumn{1}{|l|}{2.8e-2} & \multicolumn{1}{|l|}{5.1e-10} & 64 &
\multicolumn{1}{|l|}{2.4e-2} & \multicolumn{1}{|l|}{1.7e-14} & 120 \\
LSQR\_B & \multicolumn{1}{|l|}{2.8e-2} & \multicolumn{1}{|l|}{1.5e-10} & 62
& \multicolumn{1}{|l|}{2.4e-2} & \multicolumn{1}{|l|}{2.4e-15} & 120 \\
MR2 & \multicolumn{1}{|l|}{1.6e-1} & \multicolumn{1}{|l|}{3.7e-6} & 15 &
\multicolumn{1}{|l|}{} & \multicolumn{1}{|l|}{} &  \\
MINRES & \multicolumn{1}{|l|}{1.0e-2} & \multicolumn{1}{|l|}{1.2e-11} & 64 &
\multicolumn{1}{|l|}{} & \multicolumn{1}{|l|}{} &  \\
RILEY & \multicolumn{1}{|l|}{9.6e-3} & \multicolumn{1}{|l|}{8.0e-10} & 2 &
\multicolumn{1}{|l|}{1.3e-5} & \multicolumn{1}{|l|}{1.3e-10} & 2%
\end{tabular}%
\end{center}
\caption{Results for SHAW and BAART. }
\label{TABLE2}
\end{table}

The results of Tables \ref{TABLE1} and \ref{TABLE2} are of course
encouraging, especially considering the accuracy with respect to the number
of iterations. Indeed, both RA and Riley's method require a linear system to
solve at each step, and so it is fundamental to keep the number of
iterations low. However, it is worth pointing out that, in the experiments, such
linear systems are solved with the LU or Cholesky factorization, so that
most part of the computational cost is due to the first iteration.

\begin{figure}[h]
\begin{center}
\includegraphics[width=15cm]{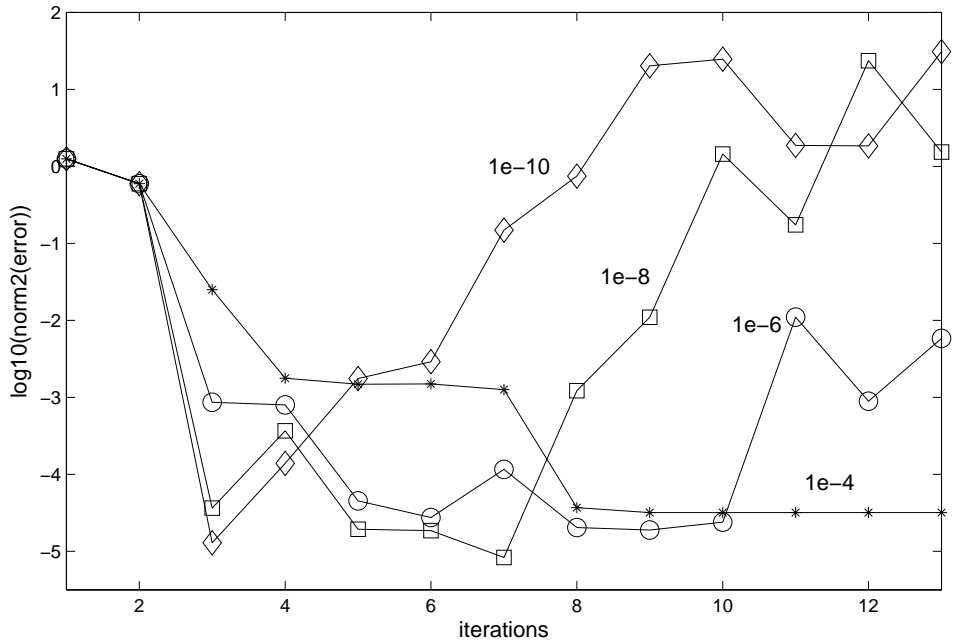}
\end{center}
\caption{BAART(120) - Error behavior for $\protect\lambda %
=10^{-4},10^{-6},10^{-8},10^{-10}$. }
\label{BBART120}
\end{figure}

A classical drawback of many iterative solvers for ill-conditioned problems
is the so-called semi-convergence (see e.g. \cite{B}), that is the
iterations initially approach the exact solution but quite rapidly diverges.
This phenomenon is very common in particular for iterative refinement methods
(thus for Riley's and RA) where there is a heavy propagation of errors. Of
course, unless a sharp error estimator is available, this undesired behavior
can be quite dangerous for applications. In order to understand what we can
do to face this problem, in Fig. \ref{BBART120} we consider the error
behavior of the RA method for BAART changing the value of the parameter.

Looking at Fig. \ref{BBART120}, we can observe that increasing $\lambda $
the procedure becomes absolutely stable, even if we have to pay a small
price in terms of accuracy. Therefore, for applications in which it is not
possible to monitor in some way the accuracy step by step, the
semi-convergence can be prevented taking $\kappa (A)^{-1/2}\ll \lambda \leq
\kappa (A)^{-1/4}$, thus looking for a compromise between accuracy and
stability. On the other side, reducing $\lambda $, the method is really fast
but also highly unstable. This last consideration is particularly true for
Riley's method, where, at least for these kind of problems, one always
observes a rapid divergence after a couple of iterations, also for
relatively large values of $\lambda $.

\begin{figure}[h]
\begin{center}
\includegraphics[width=15cm]{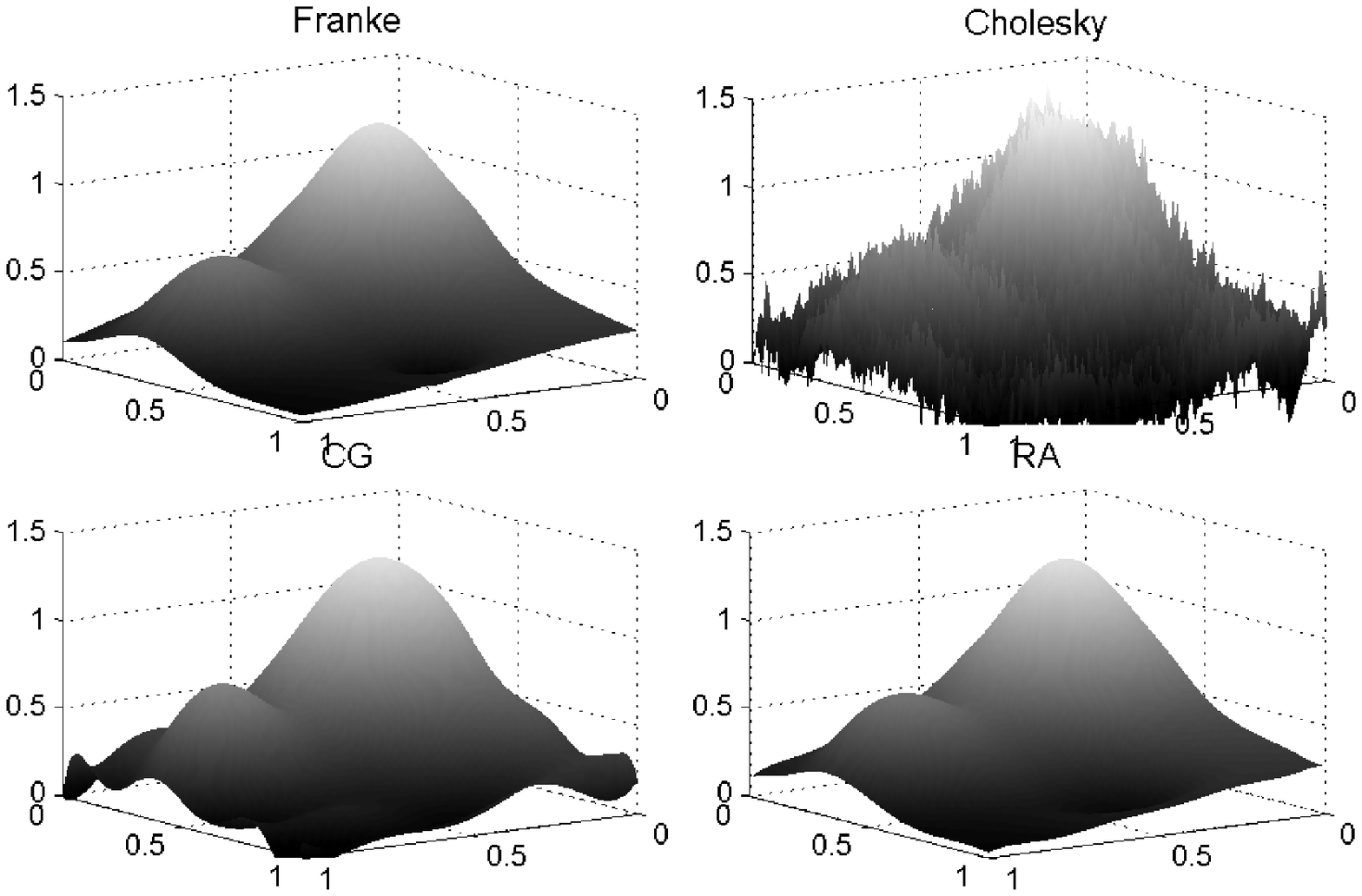}
\end{center}
\caption{Interpolation of Franke's bivariate test function by means of
Gaussian RBF. }
\label{Franke}
\end{figure}

In this Section, we also look at another classical example coming out from
approximation theory. We consider in particular the reconstruction of the
Franke's bivariate test function via interpolation by means of Gaussian
Radial Basis Functions (RBF)\ with shape coefficients equal to 1 (see e.g.
\cite{Fassa} for a background). For simplicity, instead of scattered points,
we consider here the very special case of a grid of $15\times 15$ equally
spaced points on the square $[0,1]\times \lbrack 0,1]$ that leads to a SPD
linear systems of dimension $225$ whose condition number is about $10^{21}$.
In Fig. \ref{Franke}, the surfaces obtained with the Cholesky factorization,
the CG and the RA method (with $\lambda =10^{-11}$) are plotted. Since the
exact solution of the system is unknown, we used the residual as a stopping
criterion, so that the CG result corresponds to the iteration 190 (residual $%
\approx 1.6\mathrm{e}-1$), while the RA result corresponds to the iteration
10 (residual $\approx 1.4\mathrm{e}-1$).

While the result with the Cholesky factorization was expected (a similar test
have been presented in \cite{Fass}), the difficulties with Krylov methods were
not. Indeed, the CG method has shown to be the best Krylov method for this
problem, but the results are poor if compared with those of the RA method. We
have to point out that, for this case, the reconstruction given by the RA
and the Riley's method are very similar.

\section{Extension to Tikhonov regularization}

\label{sette}

In many applications it is often necessary to deal with ill-conditioned
linear systems in which the right hand side is affected by noise. Defining $%
e_{b}$ as a perturbation (of course unknown) of the right hand side $b$, one
is forced to solve in some way%
\begin{equation}
A\widetilde{x}=\widetilde{b},\quad \widetilde{b}:=b+e_{b},  \label{pp}
\end{equation}%
hoping that the computed solution of (\ref{pp}) is close to the solution of $%
Ax=b$. In this situation, the RA method does not seem to be so powerful and
robust as in the noise-free case. Moreover, unless the noise level is very
low, it is also difficult to design a strategy to define the parameter $%
\lambda $. Indeed, in order to adopt the theory of Section \ref{cinque}
based on the analysis of the conditioning, we should need, for instance, to
construct an invertible linear filter $F$ such that $Fe_{b}\approx 0$. In
this way $F^{-1}Ax\approx \widetilde{b}$, and hence information on the choice
of $\lambda $ can be obtained considering $\kappa (F^{-1}A)$. Anyway this
kind of approach is beyond the purpose of this paper, and we prefer to extend
the idea of the RA method in order to make it able to work directly with
Tikhonov regularization in its standard form.

As well known Tikhonov regularization is based on the solution of the
minimization problem%
\begin{equation}
\min_{x}\left( \left\Vert Ax-\widetilde{b}\right\Vert ^{2}+\lambda
\left\Vert Hx\right\Vert ^{2}\right) ,\quad \lambda >0,  \label{mn}
\end{equation}%
where the matrix $H$ is generally taken as an high-pass filter (e.g. the
second derivative) so that the term $\left\Vert Hx\right\Vert ^{2}$ plays
the role of the penalization term in a constrained minimization. The main
problem is that the noise generally involves also frequencies of the exact
solution so that it is not possible to solve (\ref{mn}) letting $\lambda
\rightarrow \infty $ as in standard constrained minimization. Anyway,
defining suitably $\lambda $ (see \cite{PCH} for a background), the
corresponding solution $x_{\lambda }$ is expected to be somehow similar to
the desired noise-free solution. The problem (\ref{mn}) leads to the
solution of the regularized system%
\begin{equation}
(A^{T}A+\lambda H^{T}H)x_{\lambda }=A^{T}\widetilde{b},  \label{tik}
\end{equation}%
where the matrix $A^{T}A+\lambda H^{T}H$ is also expected to be better
conditioned than $A$.

Following the idea of the RA method, we consider here the transformation%
\begin{equation*}
Z=(A^{T}A+\lambda H^{T}H)^{-1}.
\end{equation*}%
Since the exact solution can be written as $x=\left( A^{T}A\right)
^{-1}A^{T}b$, we have%
\begin{eqnarray*}
x &=&\left( Z^{-1}-\lambda H^{T}H\right) ^{-1}A^{T}b, \\
&=&f(Q)\left( H^{T}H\right) ^{-1}A^{T}b,
\end{eqnarray*}%
where%
\begin{equation*}
Q=Z\left( H^{T}H\right) =\left( \left( H^{T}H\right) ^{-1}A^{T}A+\lambda
I\right) ^{-1}.
\end{equation*}%
Note that we are assuming to work with the exact right hand side even if, in
practice, the method is applied with $\widetilde{b}$.

Hence we can compute the solution working with the Arnoldi algorithm based
on the construction of the Krylov subspaces $K_{m}(Q,\left( H^{T}H\right)
^{-1}A^{T}b)$. Thus, starting from $v_{1}=v/\left\Vert v\right\Vert $, where $%
v$ is the solution of%
\begin{equation}
\left( H^{T}H\right) v=A^{T}b,  \label{fs}
\end{equation}%
we need to compute, at each step of the algorithm, the vectors $w_{j}=Qv_{j}$%
, $j\geq 1$, that is, we need to solve systems of the type%
\begin{equation*}
(A^{T}A+\lambda H^{T}H)w_{j}=\left( H^{T}H\right) v_{j}\text{.}
\end{equation*}%
Note that by (\ref{fs}) and the arising definition of $v_{1}$, the first
step of the Arnoldi algorithm yields the Tihhonov regularized solution $%
x_{\lambda }$ (cf. (\ref{tik})). Hence, also in this case, the procedure can
be interpreted as an iterated Tikhonov regularization.

In order to appreciate the potential of this extension (that we indicate by
RAT, Rational-Arnoldi-Tikhonov) we consider the test problem SHAW and BAART
with a right hand side contaminated by an error $e_{b}$ defined by%
\begin{equation*}
e_{b}=\frac{\delta \left\Vert b\right\Vert }{\sqrt{N}}\;u,
\end{equation*}%
where $\delta $ is the relative noise level, and $u$ is a vector containing
random values drawn from a normal distribution with mean $0$ and standard
deviation $1$. In the experiments, we define $\delta =10^{-3}$, and, as
suggested in \cite{CRS}, we take as regularization matrix%
\begin{equation*}
H=\left(
\begin{array}{ccccc}
2 & -1 &  &  &  \\
-1 & 2 & -1 &  &  \\
& \ddots & \ddots & \ddots &  \\
&  & -1 & 2 & -1 \\
&  &  & -1 & 2%
\end{array}%
\right) \in \mathbb{R}^{N\times N}.
\end{equation*}%
Indeed, at least for these experiments, this choice produces better results
than the classical $(N-2)\times N$ matrix representing the second derivative
operator. Since the noise is randomly generated, for both examples we
consider two tests, and we compare the RAT method (with different values of
the parameter $\lambda $) with GMRES, ART, LSQR\_B and MR2. The results are
collected in Table \ref{TABLE3}.

\begin{table}[th]
\begin{center}
\begin{tabular}{l|l|l|l|l|l|l|l|l|l|}
\multicolumn{2}{c}{~} & \multicolumn{4}{|c|}{SHAW(64)} & \multicolumn{4}{c|}{
BAART(120)} \\ \hline\hline
\multicolumn{1}{c}{~} & \multicolumn{1}{c|}{~} & \multicolumn{2}{c|}{test \#1
} & \multicolumn{2}{c|}{test \#2} & \multicolumn{2}{c|}{test \#1} &
\multicolumn{2}{c|}{test \#2} \\ \hline
& \multicolumn{1}{c|}{$\lambda $} & \multicolumn{1}{c|}{\textit{err}} &
\multicolumn{1}{c|}{\textit{nit}} & \multicolumn{1}{c|}{\textit{err}} &
\multicolumn{1}{c|}{\textit{nit}} & \multicolumn{1}{c|}{\textit{err}} &
\multicolumn{1}{c|}{\textit{nit}} & \multicolumn{1}{c|}{\textit{err}} &
\multicolumn{1}{c|}{\textit{nit}} \\ \hline
{\bf RAT} & 1e-3 & 0.287 & 5 & 0.215 & 3 & 0.046 & 2 & 0.046 & 2 \\
& 1e-2 & 0.293 & 5 & 0.242 & 5 & 0.028 & 3 & 0.035 & 3 \\
& 1e-1 & 0.226 & 9 & 0.230 & 7 & 0.022 & 3 & 0.029 & 3 \\
& 1e-0 & 0.297 & 7 & 0.269 & 8 & 0.010 & 3 & 0.013 & 3 \\
& 1e+1 & {\bf 0.199} & 14 & 0.269 & 8 & {\bf 0.007} & 3 & 0.009 & 3 \\
& 1e+2 & 0.293 & 18 & {\bf 0.173} & 10 & 0.008 & 4 & {\bf 0.007} & 3 \\
& 1e+3 & 0.288 & 11 & 0.268 & 13 & 0.008 & 4 & 0.010 & 4 \\
& 1e+4 & 0.575 & 10 & 0.522 & 7 & 0.008 & 4 & 0.010 & 4 \\
GMRES &  & 0.392 & 7 & 0.374 & 7 & 0.059 & 3 & 0.056 & 3 \\
ART &  & 0.837 & 64 & 0.837 & 11 & 0.344 & 120 & 0.340 & 120 \\
LSQR\_B &  & 0.361 & 14 & 0.375 & 10 & 0.142 & 6 & 0.147 & 4 \\
MR2 &  & 0.355 & 12 & 0.288 & 9 &  &  &  &
\end{tabular}%
\end{center}
\caption{Minimum attained error and corresponding iteration number for SHAW
and BAART with Gaussian noise of level $\protect\delta =10^{-3}$ }
\label{TABLE3}
\end{table}

Similarly to the noise-free case, we also consider the stabilizing effect of
a careful choice of $\lambda $. Indeed, in Figure 4 we plot the error
behavior of some of the methods considered for the solution of SHAW(64).
Taking $\lambda =10$ for the RAT method, we can overcome the problem of
semi-convergence keeping at the same time a good level of accuracy contrary
to other well performing methods such as GMRES and LSQR\_B.

\begin{figure}[h]
\begin{center}
\includegraphics[width=15cm]{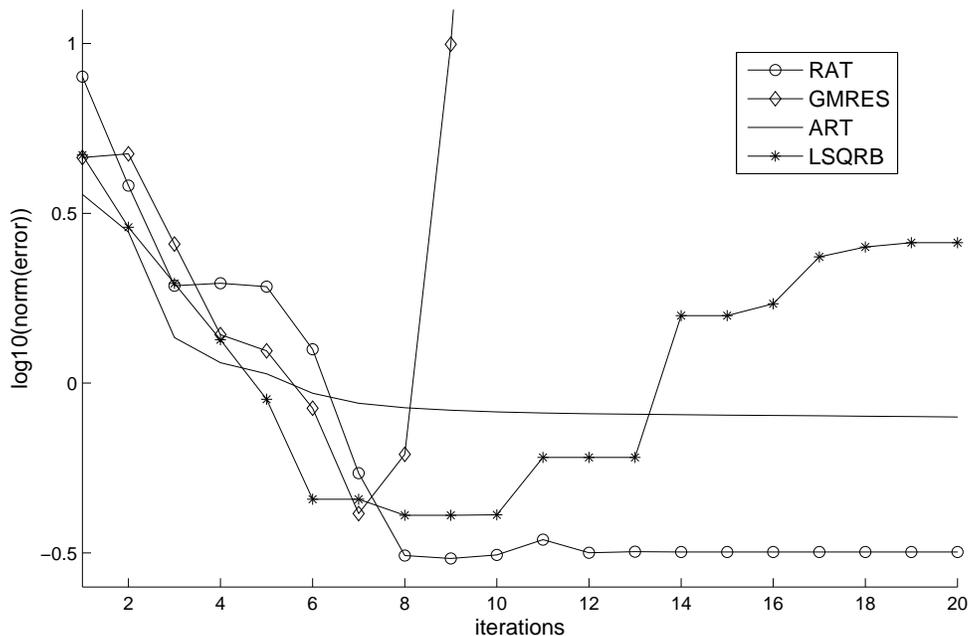}
\end{center}
\caption{Error behavior for SHAW(64) with noise. RAT method is implemented
with $\protect\lambda =10$.}
\end{figure}

\section{Conclusions}

Our experience with the RA and the RAT methods leads us to consider these
methods as reliable alternatives to the classical iterative solvers for
ill-conditioned problems. Since they actually are iterative refinement
processes, the attainable accuracy is almost never worse that the other
solvers. While this property could be somehow expected, maybe the most
important feature of these methods is their robustness. Indeed, contrary to
other iterative refinement processes such as the Riley's algorithm, the
methods work pretty well for a large window of values of $\lambda $. Hence,
having a good error estimator or working with applications in which it is
possible to monitor the result step by step, one may reduce $\lambda $ in
order to save computational work; in the opposite case, one may increase $%
\lambda $ slowing down the method but assuring a stable convergence.
To this purpose, we intend to use, in a forthcoming work, the estimates of the norm of the error described in \cite{CBerr} and \cite{CBerr2} which are based on an extrapolation procedure of the moments of the matrix of the system with respect to the residuals of the iterative method.

\vspace{0.4cm} \noindent\textbf{Acknowledgement:} The authors are grateful
to Marco Donatelli, Igor Moret, Giuseppe Rodriguez, and Marco Vianello for
many helpful discussions and comments.

\end{document}